\documentclass[
leqno,          
12pt,           
$draft,          
a4paper         
]{amsart}

\usepackage[colorlinks,citecolor=blue,urlcolor=blue]{hyperref}          
\usepackage{fourier}                                                    
\usepackage{amssymb, amsmath, amsfonts, amsthm}                         
\usepackage{mathtools, cite, enumerate, rotating, framed, lipsum}       
\usepackage{fancybox, bbm}                                              
\usepackage[dvipsnames]{xcolor}                                         
\usepackage[polish, english]{babel}                                     
\usepackage[T1]{fontenc}                                                
\usepackage[cp1250]{inputenc}                                           
\usepackage[normalem]{ulem}
\usepackage{yfonts}
\usepackage{mathrsfs}

\DeclareMathAlphabet{\mathpzc}{OT1}{pzc}{m}{it}
 \numberwithin{equation}{section}                        

\newcommand{\thmcount}{equation}                 
\newcounter{specialcounter}

\newtheorem{Thm}[\thmcount]{Theorem}
\newtheorem{Sthm}[specialcounter]{Theorem}
\newtheorem{Cor}[\thmcount]{Corollary}
\newtheorem{Lem}[\thmcount]{Lemma}
\newtheorem{Prop}[\thmcount]{Proposition}
\newtheorem{Rem}[\thmcount]{Remark}
\newtheorem{Defn}[\thmcount]{Definition}
\newtheorem{Ex}[\thmcount]{Example}
\newtheorem{Asu}[\thmcount]{Assumption}
\newtheorem{Sol}[\thmcount]{Solution}
\newtheorem*{Thmx}{Theorem}
\newtheorem*{Corx}{Corollary}
\newtheorem*{Lemx}{Lemma}
\newtheorem*{Propx}{Proposition}
\newtheorem*{Remx}{Remark}
\newtheorem*{Defnx}{Definition}
\newtheorem*{Exx}{Example}
\newtheorem*{Asux}{Assumption}
\newtheorem*{Solx}{Solution}

\newcommand \eq[1]{\begin{equation} #1 \end{equation}}
\newcommand \eqx[1]{\begin{equation*}  #1 \end{equation*}}
\newcommand \al[1]{\begin{align} #1 \end{align}}
\newcommand \alx[1]{\begin{align*}  #1 \end{align*}}
\renewcommand \sp[1]{\begin{equation} \begin{split} #1 \end{split} \end{equation}}
\newcommand \spx[1]{\begin{equation*} \begin{split} #1 \end{split} \end{equation*}}
\newcommand \en[1]{\begin{enumerate}  #1 \end{enumerate}}
\newcommand \ite[1]{\begin{itemize}  #1 \end{itemize}}
\newcommand{\thm}[2]{\begin{Thm} \label{#1} #2 \end{Thm}}
\newcommand{\sthm}[2]{\begin{Sthm} \label{#1} #2 \end{Sthm}}

\newcommand{\lem}[2]{\begin{Lem} \label{#1} #2 \end{Lem}}

\newcommand{\cor}[2]{\begin{Cor} \label{#1} #2 \end{Cor}}

\newcommand{\prop}[2]{\begin{Prop} \label{#1} #2 \end{Prop}}

\newcommand{\pr}[1]{\begin{proof} #1 \end{proof}}

\newcounter{comcount}

\renewcommand{\hline}{\vbox{\hrule width\textwidth height 1pt}\smallskip}

\renewcommand{\a}{\alpha}       \newcommand{\be}{\beta}         \newcommand{\e}{\varepsilon}
\newcommand{\w}{\omega}                 \newcommand{\de}{\delta}
        \newcommand{\la}{\lambda}
\newcommand{\vphi}{\varphi}     

 \newcommand{\Bb}{\mathbf{B}} 
\newcommand{\CC}{\mathbb{C}}

  \newcommand{\hh}{\mathcal{H}}

  \newcommand{\mm}{\mathcal{M}}
\newcommand{\NN}{\mathbb{N}}

\newcommand{\RR}{\mathbb{R}}  
  
 \newcommand{\Tt}{\mathbf{T}}

\newcommand{\ZZ}{\mathbb{Z}}  

\newcommand{\supp}{\mathrm{supp}}
\newcommand{\8}{\infty}

\renewcommand{\rm}[1]{\mathrm{#1}}
\newcommand{\wt}[1]{\widetilde{#1}}

\newcommand{\abs}[1]{\left| #1 \right|}
\newcommand{\set}[1]{\left\{ #1 \right\}}
\newcommand{\norm}[1]{\left\| #1 \right\|}

\newcommand{\eee}[1]{\left( #1 \right)}

\newcommand{\pa}{\varphi_{\a}}

\newcommand{\dtt}{\frac{dt}{t}}

\renewcommand{\Bb}{B}
\renewcommand{\hh}{H}

\title[ Sharp multiplier theorem for multidimensional Bessel operators ]{ Sharp multiplier theorem for multidimensional Bessel operators }
 \author[ Edyta Kania ]{ Edyta Kania }
 \address{
 Edyta Kania \newline
 \indent Instytut Matematyczny, Uniwersytet Wroc\l awski \newline
 \indent pl. Grunwaldzki 2/4, 50-384 Wroc\l aw, Poland }
 \email{edyta.kania@uwr.edu.pl }

 \author[ Marcin Preisner ]{ Marcin Preisner }
 \address{
 Marcin Preisner \newline
 \indent Instytut Matematyczny, Uniwersytet Wroc\l awski \newline
 \indent pl. Grunwaldzki 2/4, 50-384 Wroc\l aw, Poland }
 \email{marcin.preisner@uwr.edu.pl }

\thanks{ The research is supported by the grant No. 2017/25/B/ST1/00599 from National Science Centre (Narodowe Centrum Nauki), Poland. }

\begin{document}
\begin{abstract}
Consider the multidimensional Bessel operator
$$B f(x) = -\sum_{j=1}^N \eee{\partial_j^2 f(x) +\frac{\a_j}{x_j} \partial_j f(x)}, \quad x\in(0,\8)^N.
$$
Let $d = \sum_{j=1}^N \max(1,\a_j+1)$ be the homogeneous dimension of the space $(0,\8)^N$ equipped with the measure $x_1^{\a_1}... x_N^{\a_N} dx_1...dx_N$. In the general case $\a_1,...,\a_N >-1$ we prove multiplier theorems for spectral multipliers $m(B)$ on $L^{1,\8}$ and the Hardy space $H^1$. We assume that $m$ satisfies the classical H\"ormander condition
$$\sup_{t>0} \norm{\eta(\cdot) m(t\cdot)}_{W^{2,\be}(\RR)}<\8$$
with $\be > d/2$.

Furthermore, we investigate imaginary powers $B^{ib}$, $b\in \RR$, and prove some lower estimates on $L^{1,\8}$ and $L^p$, $1<p<2$. As a consequence, we deduce that our multiplier theorem is sharp.
\end{abstract}


\maketitle


\section{Introduction and main results}\label{intro}

\subsection{The Bessel operator}

Let $N\in \NN$ and $\a = (\a_1, ... , \a_N)$, where $\a_j > -1$ for $j=1,...,N$. Consider the space $X= (0,\8)^N$ equipped with the Euclidean metric and the measure $d\nu(x) = x^\a dx = x_1^{\a_1}...x_N^{\a_N}\, dx_1 ... dx_N$.  It is well-known that $X$ satisfies the doubling property, i.e.
\eqx{
	\nu(B(x,2r)) \leq C \nu(B(x,r)), \qquad x\in X, \, r>0,
}
where $B(x,r) = \set{y\in X \ : \  |x-y|<r}$ . 
In other words,  there exist $d, C_d>0$ such that
\eq{\tag{D} \label{doublingX}
	\nu(B(x,\gamma r)) \leq C_d (1+\gamma)^d \nu(B(x,r)), \qquad x\in X, \, r, \gamma>0.
}
We choose the constant $d$ (''homogeneous dimension'') as small as possible. In this case
   \eq{\label{d}
     d = \sum_{j=1}^N \max(1,\a_j+1).
    }

The multidimensional Bessel operator is given by $B = B_1 + ... + B_N$, where
\eqx{
B_j f(x) = - \partial_j^2 f(x) - \frac{\a_j}{x_j} \partial_j f(x), \qquad x\in X.
}
The operator $B$, initially defined on, say, $(C_c^2((0,\infty)))^n$, extends to a self-adjoint operator on $L^2(X)$. Slightly abusing notation, we shall denote this extension by the same symbol $B$. For a precise definition of $B$ we refer the reader to e.g.\cite[Sec. 2]{BCN} (see also \cite{Muckenhoupt_Stein}). Also, $B$ is the infinitesimal generator of the Bessel semigroup $\Tt_t f(x) = \int_X T_t(x,y) f(y)\, d\nu(y)$, where $T_t(x,y) = T_t^{[1]}(x_1,y_1) \cdot ... \cdot T_t^{[N]}(x_N,y_N)$ and
\eq{\label{Kt}
T_t^{[j]}(x_j,y_j) = \frac{1}{2t} (x_jy_j)^{-(\a_j-1)/2} I_{(\a_j-1)/2}\left(\frac{x_jy_j}{2t}\right) \exp\left(-\frac{x_j^2+y_j^2}{4t}\right), \quad x_j,y_j, t>0.
}
Here $I_\tau(x) = \sum_{m=0}^\8 \frac{1}{m!\Gamma(m+\tau+1)} \left(\frac{x}{2}\right)^{2m+\tau}$ is the modified Bessel function of the first kind. The kernel $T_t(x,y)$ satisfies the upper and lower gaussian bounds, i.e. there exist constants $c_1,c_2, C_1, C_2 > 0$, such that
\eq{\tag{G}\label{gauss_bess}
C_1 \nu(B(x,\sqrt{t}))^{-1} \exp\left( - \frac{|x-y|^2}{c_1 t}\right) \leq T_t(x,y) \leq C_2 \nu(B(x,\sqrt{t}))^{-1} \exp\left( - \frac{|x-y|^2}{c_2 t}\right).
}
This fact is well known and follows from the asymptotics for $\nu(B(x,\sqrt{t}))$ and $I_\tau$. For details see e.g. \cite[Lem. 4.2]{DPW_JFAA}.

Since $B$ is self-adjoint and nonnegative, for a Borel function $m: (0,\8) \to \CC$ the spectral theorem defines the operator
\eqx{
m(B) = \int_0^\8 m(\la) \, dE_B(\la),
}
where $E_B$ is the spectral resolution of $B$.

\subsection{Multiplier theorems for $B$}
Multiplier theorems for $B$ and other operators are one of the main topics in harmonic analysis. Many authors investigated assumptions on $m$ that guarantee boundedness of $m(B)$ on various function spaces, such as $L^p(X)$, $H^p(X)$, $L^{p,q}(X)$ and others.

For example, in \cite{Gosselin_Stempak} the authors proved weak type (1,1) estimates on $m(B)$ assuming $N=1$, $\a >0$ and
	\eqx{
    \eee{\int_{R/2}^R |m^{(s)}(\la)|^2 d\nu(\la)}^{1/2} \leq C R^{(\a+1)/2-s}, \quad R>0,
    }
where $s=0,...,K$ and $K$ is the least even integer greater than $(\a+1)/2 = d/2$ (see also \cite{Kapelko}). In \cite{DP_Monats}, assuming still $N=1$ and $\a>0$, it is proved that if
\eq{\tag{S} \label{Sob}
\sup_{t>0} \norm{\eta(\cdot)m(t \cdot )}_{W^{2,\beta}(\RR)} < \8
}
with some $\be > d/2$, then $m(B)$ is bounded on the Hardy space $H^1(B)$ related to $B$. Here and thereafter $W^{2,\be}(\RR)$ is the $L^2$-Sobolev space on $\RR$ and $\eta$ is a fixed nonnegative smooth cut-off function such that $\mathrm{supp} \, \eta \subseteq (2^{-1},2)$.

In the multidimensional case $N\geq 1$ in \cite{BCC} the authors prove weak type $(1,1)$ estimates for $m(B)$, where $m$ is of Laplace transform type, i.e. there exists $\phi \in L^\8 (0,\8)$, such that
\eqx{
m(x) = |x|^2 \int_0^\8 e^{-t|x|^2} \phi(t) \, dt, \quad x\in (0,\8)^N.
}
Notice, that if $m$ is of Laplace transform type, then $m$ is radial and (as a function on $(0,\8)$) satisfies \eqref{Sob} with any $\be >0$.  Another multidimensional result can be found in \cite{DPW_JFAA}, where it is proved that $m(B)$ is weak type $(1,1)$ and bounded on the Hardy space $H^1(X)$ provided that $\a_j>1$ for $j=1,...,N$ and $m$ satisfies \eqref{Sob} with $\be >d/2$. See also e.g. \cite{Garrigos_Seeger, Gasper_Trebels, Wrobel_Hankel} for other multiplier results for the Bessel operator.

Our first main goal is to obtain multiplier theorem for $B$ in the most general case $N\geq 1$ and $\a_j >-1$, $j=1,...,N$. Let us notice that many of the results before assumed that $\a_j>0$ and the case $\a_j<0$ is more difficult and less known. One reason for that is the singularity at zero of the measure $x^{\a_j} dx_j$ when $\a_j<0$. Also, so-called ''generalized translation'' operators and convolution structure for $B$ (see, e.g. \cite[Sec. 2]{BDT_d'Analyse}), does not help when $\a_j<0$. This is strictly related to the fact, that the generalized eigenfunctions of $B$ are no longer bounded if $\a_j<0$ for some $j$ and, therefore, the generalized translation is not even bounded on $L^2$. Let us also notice that, we are interested in multiplier results that are sharp in the sense that we assume \eqref{Sob} with $\be$ as small as possible. In this case this is expected to be $\be >d/2$ (we shall discuss this in Subsection \ref{impowers} below).

To state the multiplier result let us recall that the weak $L^1$ space is given by the semi-norm
	\eqx{
    \norm{f}_{L^{1,\8}(X)} = \sup_{\la>0} \la \nu \set{  x\in X \ : \ \abs{f(x)}>\la }, }
and the Hardy space $H^1(B)$ related to $B$ can be defined by the norm
\eqx{
\norm{f}_{H^1(B)} = \norm{\sup_{t>0} \abs{\Tt_t f}}_{L^1(X)}.
}
In the case $N=1$ the space $H^1(B)$ was studied in \cite{BDT_d'Analyse}, where $H^1(B)$ was characterized by means of atomic decompositions and the Riesz transforms. In the general case $N\geq 1$ and $\a_j>-1$, $j=1,...,N$ the atomic characterization of $H^1(B)$ can be found in \cite{Dziubanski2017} (see also \cite{DPW_JFAA}, \cite{JD_JGA}). We shall recall this characterization in Subsection \ref{sec-hardy-spaces} below.

\sthm{main}{
Let $N\geq 1$ and $\a_j >-1$ for $j=1,...,N$. Assume that $m:(0,\8) \to \CC$ satisfies \eqref{Sob} with $\be >d/2$, see \eqref{d}. Then:
\en{
\item $m(B)$ is bounded from $L^{1}(X)$ to $L^{1,\8}(X)$,
\item $m(B)$ is bounded from $H^{1}(B)$ to $H^{1}(B)$,
\item $m(B)$ is bounded from $L^{p}(X)$ to $L^{p}(X)$, $1<p<\8$.

}
}

Part {\bf 1.} of Theorem \ref{main} will be proved by using results of \cite{Sikora_JFA2}. More precisely, we shall check the assumptions of \cite[Th. 3.1]{Sikora_JFA2}. The proof of {\bf 2.} will be given in Section \ref{sec2}. In fact, in the proof we shall only use general properties of $B$, such as e.g. \eqref{doublingX}, \eqref{gauss_bess}, and \eqref{plan} below. Thus, the multiplier result in Section \ref{sec2} will be formulated in a more general context. This section can be read independently of the rest of the paper and we shall use different notation. As usual, {\bf 3.} is a consequence of either {\bf 1.} or {\bf 2.} by duality and interpolation, see e.g. \cite{Bernicot}.

\subsection{Imaginary powers of $B$} \label{impowers}
Another goal of this paper is to study the imaginary powers $B^{ib}$, $b\in\RR$, of the Bessel operator and establish lower bounds of these operators on some function spaces. We shall concentrate our attention on the dependence of the lower estimates on $b$ for large $b$. This is related with sharpness of multiplier theorems and may be of independent interest. To state these estimates let us restrict ourselves to the one-dimensional case $N=1$ ($X=(0,\8)$, $d\nu(x) = x^\a\, dx$, $\a>-1$). Motivated by the identity
    \eq{\label{ewew}
    {B}^{ib} = \Gamma(-ib)^{-1} \int_0^\8 t^{-ib} e^{-t{B}} \dtt
    }
let us define for $x\neq y$ the integral kernel
\eq{\label{Kb}
K_b(x,y) = \Gamma(-ib)^{-1} \int_0^\8 t^{-ib} {T}_t(x,y) \frac{dt}{t}.
}
Notice, that the integral in \eqref{ewew} is not absolutely convergent, thus we have to explain how the kernel $K_b(x,y)$ is related to the operators $B^{ib}$. Indeed, in Subsection \ref{sec3} we shall prove that for $f\in L^\8(X)$ with compact support we have
    \eq{\label{def-im-ker}
    B^{ib} f (x) = \int_X K_b(x,y) f(y) \, d\nu(y), \qquad x\notin \mathrm{supp} f
    }

One of our goals is to provide lower estimates for $B^{ib}$.
\sthm{lower1}{
Assume that $\a>-1$. Then there exist a constant $C>0$ and a function $f$ such that $\norm{f}_{L^1(X)}=1$ and for $|b|$ large enough we have
	\eqx{
    \norm{B^{ib}f}_{L^{1,\8}(X)} \geq C|b|^{d/2} .
    }
}

\sthm{lower2}{
Assume that $\a>0$ and $p\in(1,2)$. Then there exist $C_p>0$ and $f$ such that $\norm{f}_{L^p(X)}=1$ and for $|b|$ large enough we have
	\eqx{
	\norm{B^{ib} f}_{ L^p(X)} \geq C_p |b|^{\frac{d}{2}\frac{(2-p)}{p}}.
	}
}

\newcommand{\cl}{\chi_{loc}(x,y)}
\newcommand{\cg}{\chi_{glob}(x,y)}

The proofs of Theorems \ref{lower1} and \ref{lower2} are presented in Subsection \ref{sec3}. To prove Theorem \ref{lower1} we shall carefully analyze the kernels $K_b(x,y)$. More precisely, we prove the following lemma.
\lem{kernel}{
Assume that $\a>-1$ and $b\in \RR$. Then
	\sp{ \label{kernel_eq}
    K_b(x,y)  =  &c_1(b)\eee{{x^2+y^2}}^{-ib-(\a+1)/2}\\
    & + c_2(b) (xy)^{-\a/2} |x-y|^{-2bi-1}\chi_{\{y/2<x<2y\}}(x,y)\\
    &+ c_3(b) R_b(x,y),
    }
    where
\eqx{
c_1(b) = \frac{2^{2ib+1}}{\Gamma\eee{(\a+1)/{4}}} \frac{\Gamma\eee{ib + (\a+1)/2}}{\Gamma(-ib)}, \quad c_2(b) = \frac{ 2^{2ib}}{\sqrt{\pi}} \frac{\Gamma\eee{ib+1/2}}{\Gamma\eee{-ib}}, \quad c_3(b) = \Gamma(-ib)^{-1}.
}
Moreover, there exists $C>0$ that does not depend on $b$, such that
	\eqx{
    |R_b(x,y)| \leq C xy(x+y)^{-\a-3}.
    }
}

Notice that the kernel $R_b(x,y)$ is related to an operator that is bounded on every $L^p(X)$, $1\leq p \leq \8$, uniformly in $b\in \RR$. Thus we may think of $R_b(x,y)$ as of some kind of ''error term''. However, for $|b| > 1$ the size of the constants are the following:
\eq{\label{gamma_constants}
|c_1(b)| \simeq |b|^{(\a+1)/2}, \quad |c_2(b)| \simeq |b|^{1/2}, \quad |c_3(b)| \simeq |b|^{1/2} \exp\eee{\frac{\pi |b|}{2}},
}
c.f. Lemma \ref{gamma}. Thus, $c_3(b)$ grows exponentially when $|b|\to \infty$, while the constants $c_1(b)$ and $c_2(b)$ are much smaller. It appears that the growth of the constant $c_3(b)$ will lead to a~problem in deriving lower estimates for $B^{ib}$ (since our goal is to find the exact dependence on $b$). However, we can overcome this difficulty when analyzing weak $(1,1)$ norm as in Theorem \ref{lower1}. The same trick seems not to work in other function spaces (such as $H^1(B)$, $L^p(X)$ and $L^{p,\8}(X)$ with $p>1$), thus the proof of Theorem \ref{lower2} is different and uses the integral representation of the Bessel function $I_\tau$ instead of Lemma \ref{kernel}.

As a corollary of Theorems \ref{lower1} and \ref{lower2} we obtain that Theorem \ref{main} is sharp (at least for $N=1$) in the sense that $d/2$ cannot be replaced by a smaller number. The argument is standard, but we shall present it now for the convenience of the reader. One can check that for $m_b(\la) = \la^{ib}$ we have
	\eqx{
    M_b :=\sup_{t>0} \norm{\eta(\cdot) m_b(t\cdot)}_{W^{2,\be}(\RR)} \leq |b|^\be.
    }
Also, Theorem \ref{main} actually gives that $\norm{m_b(B)f}_{L^{1,\8}(X)} \leq C M_b \norm{f}_{L^1(X)}, $ where $C$ does not depend on $b$. Combining these estimates with Theorem \ref{lower1} for $|b|$ large enough we have
	\eqx{
    {|b|^{d/2}} \leq C \norm{m_b(B)}_{L^1(X) \to L^{1,\8}(X)} \leq C {|b|^{\be}}.
    }
Therefore $\be \geq d/2$. Actually, one expects that $\be \neq d/2$, but this question is beyond the scope of this paper.

Similarly, the constant $d/2$ cannot be improved for the Hardy spaces. If $\a<0$ then $d/2=1/2$ and \eqref{Sob} with $\be<1/2$ would not even guarantee that $m$ is bounded. On the other hand, for $\a>0$ if we could prove multiplier theorem on $H^1(B)$ with a constant lower than $d/2$, then by interpolation we would have better upper bounds for $m_b(B)$ on $L^p(X)$ for $1<p<2$, which contradicts Theorem \ref{lower2} by an argument similar to the one above.

\subsection{Organization of the paper and notation.} In Section \ref{sec2} we state and prove a ,,sharp'' multiplier theorem on Hardy spaces for self-adjoint operators on spaces of homogeneous type with certain assumptions (Theorem \ref{multi2}). This is a slight generalization of Theorem \ref{main}~{\bf2.} in the spirit of \cite[Th. 3.1]{Sikora_JFA2}. In Section \ref{sec2} we shall use different notation, so that it can be read independently of the rest of the paper. In Section \ref{sec_Bess} we prove the results stated above. More precisely, first we check that $B$ satisfies assumption $(P_2)$ (see Section \ref{sec2} below) in the full generality $N\geq 1$, $\a_j >-1$ for $j=1,...,N$. Thus Theorem \ref{multi2} can be applied for $B$. Then we prove Lemma \ref{kernel} and Theorems \ref{lower1} and \ref{lower2}. We shall use standard notations, i.e. $C$ and $c$ denote positive constants that may change from line to line.

\section{Sharp multiplier theorem on Hardy spaces}\label{sec2}
\subsection{Background and general assumptions}\label{ss11}

In this section we consider a space $Y$ with a~metric $\rho$ and a nonnegative measure $\mu$. We shall assume that the triple $(Y,\rho, \mu)$ is a space of homogeneous type, i.e. there exists $C>0$ such that
$
	\mu(B(x,2r)) \leq C \mu(B(x,r)),
$
for all $x\in Y$ and $r>0$, where $B(x,r) = \set{y\in Y \ : \ \rho(x,y)<r}$, c.f. \cite{CoifmanWeiss_BullAMS}. It is well-known that this implies the existence of $d,C_d>0$ such that
\eq{\tag{D} \label{doubling}
	\mu(B(x,\gamma r)) \leq C_d (1+\gamma)^d \mu(B(x,r)), \qquad  x\in Y, \, \, r, \gamma > 0.
}
As usual, we choose $d$ as small as possible, even at the cost of enlarging $C_d$.

Let $A$ denote a self-adjoint positive operator and let $E_A$ be its spectral measure, i.e. $A = \int_0^\8 \la \, dE_A(\la)$. Denote by $\mathbf{P}_t= \exp(-tA)$ the semigroup generated by $A$. Assume that there exists an integral kernel $P_t(x,y)$ such that $\mathbf{P}_t f(x) = \int_Y P_t(x,y) f(y) \, d\mu(y)$ and that satisfies the upper gaussian bounds, i.e. there exist $c_2, C_2 > 0$ such that
\eq{\tag{UG}\label{gauss_up}
P_t(x,y) \leq C_2 \mu(B(x,\sqrt{t}))^{-1} \exp\left( - \frac{\rho(x,y)^2}{c_2 t}\right), \quad t>0, \, x,y\in Y.
}

\subsection{Multiplier theorems}  By the spectral theorem, for a Borel function $m$ on $(0,\8)$, we have the operator
\eqx{
m(A) = \int_0^\8 m(\la) \, dE_A(\la).
}

In the classical case $A=-\Delta$, $Y=\RR^D$, the H\"ormader multiplier theorem states that if $m$ satisfies \eqref{Sob} with $\be>D/2$, then $m(-\Delta)$ is weak type $(1,1)$ and bounded on $L^p(\RR^D)$ for $1<p<\8$. It is well-known that the constant $D/2$ is sharp in the sense that it cannot be replaced by a smaller constant, see e.g. \cite{Sikora_Wright}.

At this point let us recall one of many multiplier theorems on spaces of homogeneous type. Suppose $Y$ and $A$ are as in Subsection \ref{ss11}. Following \cite{Sikora_JFA2}  we introduce additional assumption. Suppose that there exists $C>0$ and $q\in [2,\8]$, such that for $R>0$ and every Borel function $m$ on $\RR$ satisfying $\supp m\subseteq [R/2, 2R]$ we have
\eq{\tag{$P_q$}\label{plan}
\int_Y \abs{K_{m(A)}(x,y)}^2 d\mu(x) \leq C \mu\left(B\left(y,R^{-1/2}\right)\right)^{-1} \norm{m(R \cdot)}_{L^q(\RR)}^2.
}
\thm{multi1}{\cite[Thm. 3.1]{Sikora_JFA2}
Assume that on a space of homogeneous type $(Y,\rho, \mu)$ there is a self-adjoint positive operator $A$ that satisfies \eqref{gauss_up}. Moreover, assume that \eqref{plan} holds with some $q\in[2,\8]$ and $m$ satisfies
\eq{\tag{$S_q$} \label{Sq}
\sup_{t>0} \norm{\eta(\cdot)m(t \cdot )}_{W^{q,\beta}(\RR)} < \8
}
with some $\be>d/2$. Then $m(A)$ is of weak type $(1,1)$ and bounded on $L^p(Y)$ for $p\in (1,\8)$.
}
At this point let us make a few comments.
\en{
\item Assuming \eqref{gauss_up} the operators $m(A)$ appearing in \eqref{plan} always have integral kernels $K_{m(A)}(x,y)$, c.f. \cite[Lem. 2.2]{Sikora_JFA2}.
\item For the Bessel operator we are interested in \eqref{Sq} and \eqref{plan} for $q=2$ only, $(S) = (S_2)$. However, in Section $2$ the results are stated and proved with an arbitrary $q\in[2,\8]$.
\item The assumption \eqref{plan} in some sense plays a role of Plancherel theorem in the proof of Theorem \ref{multi1}. It is a key to obtain the sharp range $\be >d/2$. For example, if we would allow $m$ to satisfy \eqref{Sq} with $\be >d/2+1/2$, then \eqref{plan} would be superfluous.
\item The assumption \eqref{plan} is written in \cite{Sikora_JFA2} for $m$ having support in $[0,R]$ not in $[R/2,2R]$. However, a simple inspection of the proof shows that \eqref{plan} is needed only for $m$ with $\mathrm{supp} \, m \subseteq [R/2,2R]$. This makes no difference for many operators. However, it matters e.g. when considering the Bessel operator with negative parameters $\a_j$. \label{page2}
\item Assumption \eqref{plan} in \cite{Sikora_JFA2} is written for $m(\sqrt{A})$, but we use equivalent version with $m(A)$ (therefore we replace $B(y,R^{-1})$ by $B(y,R^{-1/2})$).
}

One of the main goals of this paper is to establish a multiplier theorem on Hardy spaces. We shall use the definition of the Hardy space $H^1(A)$ associated with $A$ by means of the maximal operator of the semigroup $\mathbf{P}_t$, namely
\eqx{
H^1(A) = \set{f\in L^1(Y) \ : \ \norm{f}_{H^1(A)}:= \norm{\sup_{t>0} \abs{\mathbf{P}_t f}}_{L^1(Y)} < \8}.
}
To state our result we shall assume additionally that $P_t(x,y)$ satisfies also the the lower Gaussian bounds, namely there exist $c_1, C_1 >0$, such that
\eq{\tag{LG}\label{gauss_down}
P_t(x,y) \geq C_1 \mu(B(x,\sqrt{t}))^{-1} \exp\left( - \frac{\rho(x,y)^2}{c_1 t}\right), \quad t>0, \, x,y\in Y,
}
and that the space $(Y,\rho,\mu)$ satisfies the following assumption:
\eq{\tag{Y}\label{X}
\text{ for all } x\in Y \text{ the function } r \mapsto \mu(B(x,r)) \text{ is a bijection on } (0,\8).
}
Notice that \eqref{X} implies that $\mu(Y)=\8$ and that $\mu$ is non-atomic. Now we are ready to state the theorem.

\sthm{multi2}{
Assume that $(Y,\rho,\mu)$ is a space of homogeneous type, $d$ is as in \eqref{doubling}, and \eqref{X} is satisfied. Suppose that there is a self-adjoint positive operator $A$ such that \eqref{gauss_up}, \eqref{gauss_down}, and \eqref{plan} hold  with some $q\in[2,\8]$. If $m$ satisfies \eqref{Sq} and $\be>d/2$, then $m(A)$ is bounded from $H^1(A)$ to $H^1(A)$, i.e. there exists $C>0$, such that
	\eqx{
    \norm{m(A) f}_{H^1(A)} \leq C \norm{f}_{H^1(A)}.
    }
}

The history of multiplier theorems for spaces of homogeneous type is long and wide. The interested reader is referred to \cite{Sikora_JFA2,Carbonaro_Drag,Stein,Alexopoulos,Hormander,DPW_JFAA,Garrigos_Seeger,Sikora_et_al_JAnalMath,Christ_Trans,Meda_general,Muller_Stein,Hebisch_multipliers,Chen,DP_Argentina} and references therein. {Let us concentrate for a moment on the range of parameters $\be$ in Theorem \ref{multi2}. Obviously, in general, the range $\be>d/2$ is optimal. However, it may happen that for some particular operators one may obtain multiplier results assuming that $\be >\wt{d}/2$ with $\wt{d}<d$, see e.g. \cite{Martini_Muller,Martini_Sikora,Muller_Stein}. On the other hand, there are known families of operators for which the constant $d/2$ cannot be lower. One of the methods to prove this is to derive lower estimates for $A^{ib}$ in terms of $b\in \RR$, see \cite{Chen_Sikora,Martini_Sikora,Sikora_Wright,Sikora_JFA2}. Lastly, let us mention that some multiplier results hold also in the non-doubling case, see e.g. \cite{Cowling_AnnM83}.

Boundedness of operators on the Hardy space $H^1$ is a natural counterpart of weak type $(1,1)$ bound. For example, it is a good end point for the interpolation, see e.g. \cite{Bernicot}. However, the Hardy spaces are strictly related to some cancellation conditions and it is usually more involving to study properties of operators on the Hardy space, than on $L^p$ or $L^{p,\8}$ spaces. Let us also mention that boundedness from $H^1$ to $H^1$ obviously implies boundedness from $H^1$ to $L^1$, which is usually much easier to prove.

\subsection{Hardy spaces} \label{sec-hardy-spaces}
The Hardy spaces on spaces of homogeneous type are studied extensively from the 60's, see e.g. \cite{CoifmanWeiss_BullAMS}. In particular, now we have many atomic decompositions for $H^p$ on various spaces and operators acting on this spaces. We refer the reader to e.g. \cite{Hofmann_Memoirs,BDT_d'Analyse,Dziubanski2017,Auscher_McIntosh_Russ_JGA08} and references therein.

In this subsection we recall some results on Hardy spaces related to $A$, assuming that \eqref{doubling}, \eqref{gauss_up}, \eqref{gauss_down} and \eqref{X} are satisfied.
For the proofs and more details we refer the reader to \cite{Dziubanski2017}. Firstly, there exists the unique (up to a multiplicative constant) $A$-harmonic function $\w : Y \to \RR$ such that
    \eqx{
    C^{-1} \leq \w(x) \leq C {, \qquad x\in Y}.
    }
The function $\w$ plays a special role in the analysis of $A$ and $\mathbf{P}_t$. In particular we have the following H\"older-type estimate.
\thm{conlip}{
Suppose that the semigroup $\mathbf{P}_t$ satisfies \eqref{gauss_up}, \eqref{gauss_down}. Then there exist positive constants $\gamma, c, C$, such that if $\rho(y,z)\leq \sqrt{t}$, then
\eqx{
	\abs{\frac{P_t(x,y)}{\w(y)} - \frac{P_t(x,z)}{\w(z)}} \leq C \mu(B(x,\sqrt{t}))^{-1}\left(  \frac{\rho(y,z)}{\sqrt{t}}\right)^\gamma \exp\left( - \frac{\rho(x,y)^2}{ct} \right).}}
\noindent Theorem \ref{conlip} is quite well-known and follows from a general theory. For a short and independent proof see \cite[Sec. 4]{Dziubanski2017}.

\cor{lippt}{ There exist $\gamma, C >0$ such that if $\rho(y,z)\leq \sqrt{t}$, then
	\eqx{
	\int_Y \abs{\frac{P_t(x,y)}{\w(y)} - \frac{P_t(x,z)}{\w(z)}} \, d\mu(x) \leq C \left(  \frac{\rho(y,z)}{\sqrt{t}}\right)^\gamma.
	}
}
Using Theorem \ref{conlip} the authors of \cite{Dziubanski2017} obtained the following atomic decomposition for the elements of $H^1(A)$. Let us call a function $a: Y \to \CC$ an $(\mu,\w)$-atom, if there exists a ball $B$ in~$Y$, such that:
    \eqx{
    \mathrm{supp} \,a \subseteq B, \qquad \norm{a}_\8 \leq \mu(B)^{-1}, \qquad \int_B a(x) \w(x) \, d\mu(x)=0.
    }

\thm{hardy_atomic}{\cite[Thm. 1]{Dziubanski2017}
There exists a constant $C>0$ such that for each $f\in H^1(A)$ there exist $\la_k\in \CC$ and $(\mu,\w)$-atoms $a_k$ ($k\in \NN$), such that
    \eqx{
    f(x) = \sum_{k\in \NN} \la_k a_k(x), \quad \text{and} \quad C^{-1} \norm{f}_{H^1(A)} \leq \sum_{k\in \NN} |\la_k| \leq C \norm{f}_{H^1(A)}.
    }
}


Let us start by recalling a few consequences of \eqref{doubling} and \eqref{gauss_up}.

\lem{s21}{\cite[Lem. 2.1]{Sikora_JFA2} Suppose that \eqref{doubling} and \eqref{gauss_up} hold. Then
\eqx{
    \int_{ B(y,r)^c} |P_t(x,y)|^2 \, d\mu(x) \leq C \mu(B(y,\sqrt{t}))^{-1} \exp\left( - \frac{r^2}{c_2t}\right).
}
In particular
\eqx{ \label{norm2pt}
    \norm{P_t(x,\cdot)}^2_{L^2(Y)} \leq C \mu(B(x,\sqrt{t}))^{-1}.
}}
\lem{s41}{ \cite[Lem. 4.1]{Sikora_JFA2} For $\kappa\geq 0$ there exists a constant $C=C(\kappa)>0$ such that
    \eqx{
    \int_Y |P_{(1+i\tau)R^{-1}}(x,y)|^2 (1+ R^{1/2}\rho(x,y))^\kappa \, d\mu(x) \leq C \mu\left(B\left(y,R^{-1/2}\right)\right)^{-1} (1+|\tau|)^\kappa.
    }
}
\lem{s44}{ \cite[Lem. 4.4]{Sikora_JFA2} Suppose that \eqref{doubling} holds and $\delta > 0$. Then
\eqx{
    \int_{B(y,r)^c} (1+R^{1/2}\rho(x,y))^{-d-2\delta} \, d\mu(x) \leq C \mu\left(B\left(y,R^{-1/2}\right)\right)(1+rR^{1/2})^{-2\delta}.
}
}


\newcommand{\ol}[1]{\wt{m}}

\subsection{Key kernel estimates}
This subsection is devoted to obtain key estimates needed for the proof of Theorem \ref{multi2}. We shall assume (temporarily) that $m$ satisfies $\mathrm{supp}\, m \subseteq [R/2, 2R]$ with some $R>0$. Later we shall use a partition of unity for general $m$. Denote $m_R(\la) = m(R\la)$, so that $\mathrm{supp}\, m_R \subseteq [2^{-1},2]$. Let us notice that below the letter $q\in[2,\8]$ is always the exponent related to \eqref{plan} and \eqref{Sq}. Moreover, all the spectral operators below admit related integral kernels, which can be seen by using an argument identical as in \cite[Lem. 2.2]{Sikora_JFA2}. Let us denote $\ol{m}_t(\la) = \exp(-t\la) m(\la)$ and let $M_t(x,y)$ be the kernel associated with $\ol{m}_t(A) = \mathbf{P}_t m(A)$.

\prop{prop:main}{
	Assume that $\mathrm{supp} \, m \subseteq [R/2, 2R]$ and $m_R \in W^{q,\be}(\RR)$ with $\be > d/2$. Then, there exist $\delta, \gamma, C>0$ such that for $y,z\in Y$ and $r>0$ we have
	\eq{ \label{main1}
    \int_{B(y,r)^c} \sup_{t>0} \, \abs{M_t(x,y)} \, d\mu(x) \leq C \left(1 + rR^{1/2}\right)^{-\delta} \norm{m_R}_{W^{q,\be}(\RR)},}
    and, for $\rho(y,z) < R^{-1/2}$,
	\eq{\label{main2}
   \int_{B(y,r)^c} \sup_{t>0}\abs{ \frac{M_t(x,y)}{\w(y)} - \frac{M_t(x,z)}{\w(z)}} \, d\mu(x) \leq C \left( R^{1/2}\rho(y,z)\right)^{\gamma} \norm{m_R}_{W^{q,\be}(\RR)}.
	}
}

Let us start by showing the following lemma.

\lem{s43}{
    For $\e >0,\kappa \geq 0$ there exists a constant $C=C(\kappa,\e)$ such that
    \eqx{
    \int_Y \sup_{t>0} |M_{t}(x,y)|^2 \left(1+R^{1/2}\rho(x,y)\right)^{\kappa} \, d\mu(x) \leq C \mu\left(B\left(y,R^{-1/2}\right)\right)^{-1} \norm{m_R}^2_{W^{q,\kappa/2+\e}(\RR)}.
    }}

\pr{
Fix a cut-off function $\psi\in C_c^\8(4^{-1},4)$, such that $\psi \equiv 1$ on $[2^{-1}, 2]$. Set
\spx{
n_{t,R}(\la) = m_R(\la) \underbrace{e^{-tR\la}  e^\la \psi(\la)}_{\la_{t,R}(\la)}.
}
By the Fourier inversion formula,
    \eqx{
    \ol{m}_{t}(A) = n_{t,R}(AR^{-1})e^{-AR^{-1}} = \frac{1}{2\pi} \int_\RR \widehat{n}_{t,R}(\tau) \exp \left( (i\tau-1)AR^{-1} \right)  \, d\tau
    }
and
    \eq{\label{eq_Mt}
    M_t(x,y) = \frac{1}{2\pi} \int_\RR \widehat{n}_{t,R}(\tau) P_{(1-i\tau)R^{-1}}(x,y) \, d\tau.
    }

Notice that $\mathrm{supp} \, \la_{t,R}^{(N)} \subseteq (4^{-1}, 4)$ for arbitrary $N\in\NN$. By simple calculus we can find a~constant $C_N$ such that
\spx{
    \sup_{R>0,t>0} \abs{\widehat{\la}_{t,R}(\tau)} \leq C_N(1+|\tau|)^{-N}.
}
Since $\widehat{n}_{t,R} = \widehat{m}_R \ast \widehat{\la}_{t,R}$ and $(1+|\tau|) \leq (1+|\theta|)(1+|\tau-\theta|)$, for $\kappa\geq0$ and $\e>0$ we use the Cauchy-Schwarz inequality, getting

\spx{
\int_{\RR} \sup_{t>0} |\widehat{n}_{t,R}(\tau)|  (1+|\tau|)^{\kappa/2} \, d\tau & \leq \int_{\RR} \int_{\RR} \sup_{t>0} |\widehat{m}_R(\theta)||\widehat{\la}_{t,R}(\tau-\theta)|  (1+|\tau|)^{\kappa/2}  \, d\theta  \,d\tau \\
&  \leq \int_{\RR} \int_{\RR} \sup_{t>0} |\widehat{m}_R(\theta)||\widehat{\la}_{t,R}(\tau-\theta)|  (1+|\theta|)^{\kappa/2} \, (1+|\tau-\theta|)^{\kappa/2} \,d\tau \, d\theta \\
& \leq C  \int_{\RR}  |\widehat{m}_R(\theta)|  (1+|\theta|)^{(\kappa+1)/2+\e}  (1+|\theta|)^{-1/2-\e} \, d\theta \\
& \leq C \norm{m_R}_{W^{2,(1+\kappa)/2+\e}(\RR)} \left(  \int_{-\8}^\8  (1+|\theta|)^{-1-\e}  \, d\theta \right)^{1/2} \\
&  \leq C \norm{m_R}_{W^{2,(1+\kappa)/2+\e}(\RR)}.
}
Hence, by \eqref{eq_Mt}, the Minkowski inequality, and Lemma \ref{s41} we obtain
    \sp{ \label{W2}
    & \left( \int_Y \sup_{t>0} |M_t(x,y)|^2 (1+R^{1/2}\rho(x,y))^\kappa \ d\mu(x) \right)^{1/2} \\
    & \leq \int_\RR \sup_{t>0} |\widehat{n}_{t,R}(\tau)| \left( \int_Y |P_{(1-i\tau)R^{-1}}(x,y)|^2 (1+R^{1/2}\rho(x,y))^{\kappa} \, d\mu(x) \right)^{1/2} \, d\tau\\
    & \leq C  \mu\left(B\left(y,R^{-1/2}\right)\right)^{-1/2} \int_\RR  \sup_{t>0} |\widehat{n}_{t,R}(\tau)| (1+|\tau|)^{\kappa/2} \, d\tau \\
    & \leq C  \mu\left(B\left(y,R^{-1/2}\right)\right)^{-1/2}  \norm{m_R}_{W^{2,(1+\kappa)/2+\e}(\RR)}\\
        & \leq C  \mu\left(B\left(y,R^{-1/2}\right)\right)^{-1/2}  \norm{m_R}_{W^{q,(1+\kappa)/2+\e}(\RR)}.
    }
In the last inequality we have used that $\mathrm{supp}\, m_R\subseteq [2^{-1},2]$ and $q\geq 2$.

Observe that \eqref{W2} is exactly the estimate we look for, but the Sobolev parameter is higher by $1/2$ than we want. To sharpen this estimate, we make use of known interpolation method. Notice, that $M_t(x,y) = \mathbf{P}_t(K_{m(A)}(\cdot,y))(x)$. It is well-known that \eqref{gauss_up} implies boundedness on $L^2(Y)$ of the maximal operator $\mm f =\sup_{t>0} |\mathbf{P}_t f|$. A second estimate needed for an interpolation is the following
\sp{\label{L2}
\eee{\int_Y \sup_{t>0} |M_{t}(x,y)|^2 \, d\mu(x)}^{1/2} & = \norm{ \mm K_{m(A)}(\cdot,y)}_{L^2(Y)} \\
& \leq C \norm{  K_{m(A)}(\cdot,y) }_{L^2(Y)}\\
& \leq C \mu\left(B\left(y,R^{-1/2}\right)\right)^{-1/2}\norm{m_R}_{L^q(\RR)}.
}
In the last inequality we have used \eqref{plan}. Now, Lemma \ref{s43} follows by interpolating \eqref{W2} and \eqref{L2}, see e.g. proofs of \cite[Lem. 4.3(a)]{Sikora_JFA2} and \cite[Lem. 2.2]{DPW_JFAA} for details.
}

\pr{[Proof of \eqref{main1}] By the Cauchy-Schwarz inequality and Lemmas  \ref{s44} and \ref{s43},
 \spx{ \label{first}
    & \int_{B(y,r)^c} \sup_{t>0} |M_{t}(x,y)| \, d\mu(x) \\
    & \leq \left( \int_Y \sup_{t>0} |M_{t}(x,y)|^2(1+R^{1/2}\rho(x,y))^{d+2\de} \, d\mu(x) \right)^{1/2} \left( \int_{B(y,r)^c} (1+R^{1/2}\rho(x,y))^{-d-2\de} \, d\mu(x) \right)^{1/2}\\
    & \leq C \mu\left(B\left(y,R^{-1/2}\right)\right)^{-1/2} \norm{m_R}_{W^{q,d/2+\de+\e}(\RR)}\mu\left(B\left(y,R^{-1/2}\right)\right)^{1/2}(1+rR^{1/2})^{-\de} \\
    & \leq C (1+rR^{1/2})^{-\de} \norm{m_R}_{W^{q,\be}(\RR)},
    }
where $\de, \e >0$ are such that $d/2+\de+\e \leq \be$.
}

\renewcommand{\wt}{\widetilde}

Consider for a moment the operator $\mathbf{P}_tm(A) \exp(AR^{-1})$ and let $\wt{M}_{t,R}(x,y)$ be its kernel. By almost identical arguments as in the proofs of Lemma \ref{s43} and \eqref{main1}, we can show that for $\be>d/2$ we also have
	\eq{\label{est_Mtilde}
	\int_{B(y,r)^c} \sup_{t>0} |\wt{M}_{t,R}(x,y)| \, d\mu(x) \leq C \norm{m_R}_{W^{q,\be}(\RR)}.
	}

\pr{[Proof of \eqref{main2}]
Notice, that $M_{t}(x,y) = \int_Y \wt{M}_{t,R}(x,u)P_{R^{-1}}(u,y)  \, d\mu(u)$. For $\rho(y,z)<R^{-1/2}$, by Corollary \ref{lippt} and \eqref{est_Mtilde},
\spx{
	& \int_{B(y,r)^c} \sup_{t>0} \abs{\frac{M_t(x,y)}{\w(y)} -\frac{M_t(x,z)}{\w(z)}} \, d\mu(x) \\
	& = \int_{B(y,r)^c} \sup_{t>0} \abs{\int_Y \wt{M}_{t,R}(x,u)\left(  \frac{P_{R^{-1}}(u,y)}{\w(y)} -
	\frac{P_{R^{-1}}(u,z)}{\w(z)} \right)  \, d\mu(u)} \, d\mu(x) \\
	& \leq \int_Y  \abs{  \frac{P_{R^{-1}}(u,y)}{\w(y)} - \frac{P_{R^{-1}}(u,z)}{\w(z)}}  \int_{B(y,r)^c} \sup_{t>0} | \wt{M}_{t,R}(x,u)| \, d\mu(x) \, d\mu(u) \\
	& \leq C  \left(  R^{1/2}\rho(y,z)\right)^\gamma \norm{m_R}_{W^{q,\be}(\RR)}.
}

}

\subsection{Proof of Theorem \ref{multi2}}
{Theorem \ref{multi2} follows from Proposition \ref{prop:main} by a quite standard argument. We present the details for completeness and convenience of the reader. As usual, by a continuity argument, in order to prove boundedness of the operator $m(A)$ on $H^1(A)$ it is enough to show that there exists $C>0$ such that
\eqx{
\norm{m(A)a}_{H^1(A)} = \norm{\mm m(A)a}_{L^1(Y)} \leq C}
holds for every $(\mu,\w)$-atom $a$, see Theorem \ref{hardy_atomic}. Assume then that: $\supp \, a \subseteq B(y_0,r) = : B$, $\norm{a}_\8 \leq \mu(B)^{-1}$, and $\int a \, \w d\mu = 0$. As always, the analysis on $2B = B(y_0,2r)$ follows by the Cauchy-Schwarz inequality and boundedness of $\mm$ and $m(A)$ on $L^2(Y)$. More precisely,
\spx{
\norm{\mm m(A) a }_{L^1(2B)} & \leq \mu(2B)^{1/2} \norm{\mm m(A) a }_{L^2(Y)} \\
& \leq C \mu(B)^{1/2} \norm{ a (x)}_{L^2(Y)} \leq C.
}
Therefore, it is enough to prove that
\eq{\label{eq:main}
\norm{\mm m(A) a }_{L^1((2B)^c)} \leq C.
}
Let $\eta\in C_c^\8(2^{-1},2)$ be a fixed function such that
$
\sum_{j\in\ZZ}\eta(2^{-j}\la) = 1
$
for all $\la \in(0,\8)$. By using this partition of unity, we decompose $m$ as
\eqx{
	m(\la) = \sum_{j\in \ZZ} \eta(2^{-j}\la) m(\la) = \sum_{j\in\ZZ} m_j(\la).
}
Fix $N\in \ZZ$ such that
    $
    2^{-N} \leq r^2 < 2^{-N+1}.
    $
Then
    \spx{
   \norm{\mm m(A)a}_{L^1((2B)^c)} & \leq \sum_{j\in\ZZ} \norm{\mm m_j(A)a}_{L^1((2B)^c)}= \sum_{j\geq N}... + \sum_{j<N}...  = S_1 + S_2.
    }

    Denote $m_{j,t}(\la) = \exp(-t\la) m_j(\la)$ and let $M_{j,t}(x,y)$ be the kernel of $m_{j,t}(A) = \mathbf{P}_t m_j(A)$. Obviously, $\supp \, m_{j,t} \subseteq [2^{j-1}, 2^{j+1}]$ and applying \eqref{main1} we obtain that
\spx{ \label{s1}
S_1
& \leq  \sum_{j\geq N} \int_{(2B)^c} \int_{B} \sup_{t>0} |M_{j,t}(x,y)||a(y)| d\mu(y) d\mu(x) \\
& \leq  \sum_{j\geq N} \int_{B} |a(y)| \int_{B^c} \sup_{t>0}  |M_{j,t}(x,y)| d\mu(x)  d\mu(y)\\
& \leq C \norm{a}_{L^1(Y)} \sum_{j\geq N} (1+2^{j/2}r)^{-\delta} \norm{\eta(\cdot) m(2^j \cdot) }_{W^{q,\beta}(\RR)} \\
& \leq C \sup_{t>0} \norm{\eta(\cdot) m(t \cdot) }_{W^{q,\beta}(\RR)} \leq C.}
If $y\in B$ and $j < N$, then $\rho(y,y_0)<  r < 2^{-j/2}$ and we can apply \eqref{main2} for the kernel $M_{j,t}$  with $R=2^j$. Using the cancellation condition of $a$,
\spx{ \label{s2}
S_2 & \leq \sum_{j < N} \int_{(2B)^c} \sup_{t>0} \abs{\int_{B} M_{j,t}(x,y) a(y) d\mu(y)} d\mu(x) \\
& = \sum_{j < N} \int_{(2B)^c} \sup_{t>0} \abs{\int_{B} \left( \frac{M_{j,t}(x,y)}{\w(y)}- \frac{M_{j,t}(x,y_0)}{\w(y_0)}\right) a(y) \w(y) d\mu(y)} d\mu(x)\\
& \leq  \sum_{j < N} \int_{B}  |a(y)| \int_{B(y,r)^c} \sup_{t>0} \abs{ \frac{M_{j,t}(x,y)}{\w(y)}- \frac{M_{j,t}(x,y_0)}{\w(y_0)}} d\mu(x) \, \w(y) d\mu(y)\\
& \leq C \sum_{j< N} 2^{\frac{j\gamma}{2}} \int_{B} |a(y)| \rho(y,y_0)^{\gamma} d\mu(y) \norm{\eta(\cdot) m(2^j \cdot) }_{W^{q,\beta}(\RR)} \\
& \leq C  \sup_{t>0} \norm{\eta(\cdot) m(t \cdot) }_{W^{q,\beta}(\RR)} r^\gamma \sum_{j< N} 2^{\frac{j\gamma}{2}} \leq C.
}
This finishes the proof of \eqref{eq:main} and Theorem \ref{multi2}.

\section{The multidimensional Bessel operator}\label{sec_Bess}
In this Section we turn back to the analysis related to $B$ and prove the results stated in Section \ref{intro}.

\subsection{The Hankel transform}

Recall that $N\in \NN$ and $\a_j >-1$ for $j=1,...,N$. For $x,\xi\in X = (0,\8)^N$ denote $\pa(x\xi) = \vphi_1(x_1\xi_1)\cdot ... \cdot \vphi_N(x_N\xi_N)$, where $$\vphi_j(z) = 2^{(\a_j-1)/2}\Gamma\left((\a_j+1)/2 \right)z^{-(\a_j-1)/2}J_{(\a_j-1)/2}(z), \quad z>0.$$ Here $J_\tau$ denotes the Bessel function of the first kind. By the asymptotics of $J_\tau$ one has
	\eq{\label{phi}
   \abs{\vphi_j(z)} \leq C (1+z)^{-\a_j/2}, \quad z>0.
    }
    The Hankel transform is defined by
\eq{\label{hankel}
	H_\a f(\xi) = \int_X f(x)\pa(x\xi) d\nu(x), \quad \xi\in X,
}
As we have already mentioned, $\varphi_j \in L^\8$ if and only if $\a_j \geq 0$. Nevertheless, it is known that $H_\a$ always extends uniquely to an isometric isomorphism on $L^2(X)$, see \cite{BCC} and \cite[Lem. 2.7]{Betancor_Stempak}. The multipliers $m(B)$ and $\hh_\a$ are related in the same way, as $m(-\Delta)$ and the Fourier transform on $\RR^D$. In particular, if
\eqx{
n(\la) = m(|\la|^2), \qquad \la \in X,
}
then $m$ is radial and
\eq{\label{m-han}
m(B) f = H_\a (n \cdot H_\a).
}

\subsection{$(P_2)$ for multidimensional Bessel operator}
Let us first recall that $\mathbf{T}_t$ satisfies \eqref{gauss_bess} and, obviously, $X$ satisfies \eqref{X}. Therefore, Theorem \ref{main} follows from Theorems \ref{multi1} and \ref{multi2} provided that \eqref{plan} holds with $q=2$, which we now prove. The case $N=1$ follows by similar and simpler argument, thus we shall concentrate on $N\geq 2$. Let $k\in\set{1,...,N-1}$ and $c_j<2^{-N}$ for $j=1,...,k$. Define the sets
$$S_{c_1,...,c_k} = \set{x\in X \ : \ 1/2<|x|^2<2 \text{ and } x_{j}<c_j \text{ for } j=1,...,k}.$$

\lem{x_k}{ Suppose that $\supp \,m\subset [R/2,2R]$, $N\geq 2$, $k\leq N-1$, and $c_j <2^{-N}$ for $j=1,2,...,k$. Then there exists $C>0$ such that
\eqx{
	\int_{S_{c_1,...,c_k}} \abs{m\left(R|x|^2\right)}^2 x_{1}^{\a_{1}}  ...  x_k^{\a_k} \, dx_1...dx_N \leq C c_1^{\a_{1}+1}  ...  c_k^{\a_{k}+1} \norm{m(R\cdot)}_{L^2(\RR)}^2.
}}

\pr{
	Introduce the spherical coordinates $(r,\theta_1, ... , \theta_{N-1})$ on $\RR^N$, namely
	\eqx{ \label{spher-var-d}
	\begin{cases}
    x_1 = r\sin(\theta_1), & \\
    x_i = r \sin(\theta_i) \Pi_{j=1}^{i-1}\cos(\theta_j), \qquad \ \ \   \text{ for } i=2,3,...,N-1, \\
    x_N = r \Pi_{j=1}^{N-1}\cos(\theta_j)\\
    dx_1\,...\,dx_N = r^{N-1} \prod_{j=1}^{N-2} \cos^{N-1-j}(\theta_j) dr\, d\theta_1 ... \,d\theta_{N-1}. \end{cases}
}
Since $x\in X$, then $\theta_j \in (0,\pi/2)$ for $j=1,...,N-1$. We claim that if $x\in S_{c_1,...,c_k}$, then

\eq{\label{male_katy}
\sin(\theta_j) < 2^{-N/2} \leq 2^{-1/2}
}
for $j=1,...,k$. Obviously, if $\sin(\beta) < 2^{-1/2}$ for $\be \in (0,\pi/2)$, then $(\cos \beta)^{-1} < 2^{1/2}$. Observe that $2^{-1/2}<r<2^{1/2}$, since $\supp \, m \subseteq [R/2, 2R]$. Therefore, \eqref{male_katy} follows easily by induction, i.e. for $i=1,...,k$,
\eqx{
\sin \theta_i = x_i r^{-1} (\cos  \theta_1)^{-1} ... (\cos  \theta_{i-1})^{-1}\leq 2^{-N} 2^{i/2} \leq 2^{-N/2}.
}

 Denote $S=S_{c_1,...,c_k}$. As a consequence of \eqref{male_katy} we have that $\sin \theta_i \simeq \theta_i$ and $\cos \theta_i \simeq C$ for $i=1,...,k$. Using this, $x_i^{\a_i} \simeq r^{\a_i} \theta_i^{\a_i}$ for $i=1,...,k$ and
\spx{
	& \int_{S} \abs{m(R|x|^2)}^2 x_{1}^{\a_{1}} ... x_k^{\a_k} \, dx_1...dx_N \\
	& \leq C \int_S \abs{m(Rr^2)}^2 r^{N-1+\a_1+...+\a_k} \theta_1^{\a_1}...\theta_k^{\a_k} dr d\theta_{1} ... d\theta_{N-1} \\
	& \leq C\int_0^{c c_1} \theta_1^{\a_1} \,d\theta_1 \cdot ... \cdot \int_0^{c c_k} \theta_k^{\a_k} \,d\theta_k \cdot  \int_{1/2<r^2<2} \abs{m\left(Rr^2\right)}^2 dr \\
	& \leq C c_1^{\a_1+1} \cdot ... \cdot c_k^{\a_k+1} \cdot \norm{m(R\cdot)}_{L^2(\RR)}^2.
}}

\prop{prop_P}{
Assume that $N\in \NN$ and $\a_j>-1$ for $j=1,..,N$. Then \eqref{plan} holds for $B$ {with} $q=2$.
}

\pr{\label{plan-bessel}
In the proof we consider only the case $N\geq 2$. Let $q=2$ and suppose that $m$ is supported in $[R/2,2R]$ for some $R>0$. Notice that by \eqref{m-han} and \eqref{hankel} we have
\spx{
	m(B)f(x) & = \int_{X} f(y) \int_{X} n(\xi) \pa(y\xi) \pa(x\xi) d\nu(\xi) \, d\nu(y)\\
	& = \int_{X} f(y) H_\a\left( n(\cdot) \pa(y\cdot)\right)(x) \, d\nu(y)
}
and the kernel associated with $m(B)$ has the form
\eqx{K_{m(B)}(x,y) =  H_\a\left( n(\cdot) \pa(y\cdot)\right)(x).}
Therefore, by the Plancherel identity for $H_\a$, $(P_2)$ is equivalent to
\eq{\label{upr}
	\int_{R/2<|x|^2<2R}\abs{m(|x|^2)}^2\abs{ \pa(xy)}^2 d\nu(x) \leq C \nu\left(B\left(y,R^{-1/2}\right)\right)^{-1} \norm{m(R\cdot)}^2_{L^2(\RR)}.
}
For each $i=1,...,N$ we consider four cases:
\begin{enumerate}[{\bf C1.}]
\item $x_i < 2^{-N}$,\quad   $\sqrt{R}y_i > 2^N $, \quad $\sqrt{R}y_ix_i < 1$,

\item $x_i < 2^{-N}$,\quad  $\sqrt{R}y_i \leq 2^N $,

\item $x_i < 2^{-N}$,\quad  $\sqrt{R}y_{i} > 2^N$, \quad $\sqrt{R}y_ix_i \geq 1$,

\item $ 2^{-N} \leq  x_i \leq \sqrt{2}  $.
\end{enumerate}

Divide the set $\set{x\in X \ : \  1/2<|x|^2<2}$ into several disjoint regions using the cases above. Without loss of generality we may consider the set $S$ of points $x\in X$ such that:
\ite{
\item $x_i$ satisfies {\bf C1.} for $i=1,..., k_1$,
\item $x_i$ satisfies {\bf C2.} for $i=k_1+1,..., k_2$,
\item $x_i$ satisfies {\bf C3.} for $i=k_2+1,..., k_3$,
\item $x_i$ satisfies {\bf C4.} for $i=k_3+1,..., N$,
}
where $0\leq k_1\leq k_2 \leq k_3 < N$. The fact that $k_3< N$ is implied by $|x|^2>1/2$. Notice that it may happen that $S$ is empty. Recall that $\nu(B(y,r)) \simeq \Pi_{j=1}^N \nu_j(B(y_j,r))$, where $d\nu_j(x_j) = x_j^{\a_j} dx_j$ is the one-dimensional measure, and
\eqx{
\nu_j\left(B\left(y_j, R^{-1/2}\right)\right)^{-1} \simeq R^{(\a_j+1)/2} \eee{1+\sqrt{R}y_j}^{-\a_j}.
}

Denote $d_{gl} = N+\a_1 +...+\a_N$. Using \eqref{phi} and  Lemma \ref{x_k} with $k=k_2$, we have
\spx{
&\int_{R/2<|x|^2<2R}\abs{m\left(|x|^2\right)}^2\abs{ \pa(xy)}^2 d\nu(x)
\leq C \sum_S R^{\frac{d_{gl}}{2}} \int_S \abs{m\left(R|x|^2\right)}^2 \prod_{j=1}^N\eee{x_j^{-1}+\sqrt{R}y_j}^{-\a_j} dx\\
&\leq  C \sum_S R^{\frac{d_{gl}}{2}} \int_S \abs{m\left(R|x|^2\right)}^2 x_1^{\a_1}...x_{k_2}^{\a_{k_2}} \prod_{j=k_2+1}^N \left(1+\sqrt{R}y_{j}\right)^{-\a_{j}} dx_1...dx_N\\
&\leq C \sum_S R^{\frac{d_{gl}}{2}} \prod_{i=1}^{k_1}\left(\sqrt{R}y_i\right)^{-\a_i-1} \prod_{j=k_2+1}^N \left(1+\sqrt{R}y_{j}\right)^{-\a_{j}} \norm{m(R\cdot)}_{L^2(\RR)}^2\\
&\leq C \prod_{i=1}^{k_1}R^{(\a_i+1)/2}\left(\sqrt{R}y_i\right)^{-\a_i} \prod_{k=k_1+1}^{k_2} R^{(\a_k+1)/2} \prod_{j=k_2+1}^N R^{(\a_j+1)/2} \left(1+\sqrt{R}y_{j}\right)^{-\a_{j}} \norm{m(R\cdot)}_{L^2(\RR)}^2\\
&\leq C \nu\left(B\left(y,R^{-1/2}\right)\right)^{-1} \norm{m(R\cdot)}_{L^2(\RR)}^2.
}
}

\subsection{Imaginary powers of $B$}\label{sec3} In this subsection we prove Lemma \ref{kernel} and Theorems \ref{lower1} and~\ref{lower2}. From now on we consider one-dimensional Bessel operator, i.e. $N=1$, $X=(0,\8)$, $\a>-1$, and $d\nu(x) = x^\a \, dx$.

Let us start this section by recalling well-known asymptotics of the Bessel function $I_\tau$ \cite{Watson,Lebedev}, i.e.
	\al{\label{Bessel_loc}
    &&&&I_{\tau}(x) &= \Gamma\eee{\frac{\tau+1}{2}}^{-1} \eee{\frac{x}{2}}^\tau + O\eee{x^{\tau+1}}, && x \sim 0&&&&\\
    \label{Bessel_glob}
    &&&&I_{\tau} (x) &= (2\pi x)^{-1/2} e^{x}\eee{ 1 + O(x^{-1})}, && x\sim \8&&&&
    }

Now we provide a short argument for \eqref{def-im-ker}. In \cite[Sec. 4.3]{BCN} it is proved that $B^{ib}$ is associated with the kernel $$-\Gamma(-ib+1)^{-1} \int_0^\8 t^{ib} \partial_t T_t(x,y) \, dt$$ in the sense as in \eqref{def-im-ker} (let us notice that in \cite{BCN} only positive values of $\a$'s are considered, but the proof works for $\a_j>-1$ as well). By integrating by parts,
\spx{
- \Gamma(-ib+1)^{-1} \int_0^\8 t^{-ib} \partial_t T_t(x,y)\,dt  =& - \Gamma(-ib+1)^{-1} \lim_{\e\to 0}  \left( \e^{ib} T_{\e^{-1}}(x,y) - \e^{-ib} T_\e(x,y) \right)\\
& +  \Gamma(-ib)^{-1}   \int_0^{\8} t^{-ib} T_t(x,y)\,\dtt = K_b(x,y) .
}

\pr{[Proof of Lemma \ref{kernel}]
Let us first notice that for $\kappa\in \RR$ and $c,M>0$, there exists $C=C(\kappa,c,M)$ such that
\sp{\label{int-conv}
\int_{c z}^{\8} t^\kappa \exp\eee{-\frac{t}{4}} \dtt \leq Cz ^{-M}, \qquad z\geq 1.
}

Using \eqref{Kb} and \eqref{Kt} one obtains
	\spx{
    2 \Gamma(-ib) K_b(x,y) & = \int_0^\8 t^{-ib-1} (xy)^{-(\a-1)/2} I_{(\a-1)/2}\left(\frac{xy}{2t}\right) \exp\left( -\frac{x^2+y^2}{4t}\right)  \dtt \\
    &= \int_0^{xy} ... +  \int_{xy}^\8  ...= A_1 + A_2.
    }
Denote $\cl=\chi_{\set{y/2<x<2y}}(x,y)$ and $\cg = 1-\cl$, $x,y\in X$. In the proof below all expressions denoted by $R_k$ shall be parts of the kernel $R_b(x,y)$. Using \eqref{Bessel_glob}, we write $A_1 = A_{1,1}+R_1$, where
\spx{
A_{1,1} = &\pi^{-1/2} \int_0^{xy} t^{-ib-1/2} (xy)^{-\a/2} \exp\left( -\frac{|x-y|^2}{4t}\right) \dtt
}
and
\alx{
     |R_1| &= \abs{ \int_0^{xy} t^{-ib-1} (xy)^{-(\a-1)/2} \exp\left( -\frac{x^2+y^2}{4t}\right) \left(I_{(\a-1)/2}\left(\frac{xy}{2t}\right) -  \left(\frac{\pi xy}{t}\right)^{-1/2}\exp\left(\frac{xy}{2t}\right)\right)   \dtt}\\
     &\leq C \int_0^{xy} t^{1/2} (xy)^{-\a/2-1} \exp\eee{-\frac{|x-y|^2}{4t}} \dtt \\
     &= C |x-y| (xy)^{-\a/2-1}  \int_{\frac{|x-y|^2}{xy}}^\8 t^{-1/2} e^{-t/4} \dtt\\
    &\leq C xy (x+y)^{-\a-3}.
}
In the last inequality we have used \eqref{int-conv}. Denoting
\eqx{
R_2 =  \cg A_{1,1} \quad \text{and}\quad R_3 = \pi^{-1/2} \cl \int_{xy}^\8 t^{-ib-1/2} (xy)^{-\a/2} \exp\left( -\frac{|x-y|^2}{4t}\right) \dtt
}
we have
\spx{
A_{1,1} - R_2 + R_3 =&\pi^{-1/2}\cl \eee{\int_0^\8 t^{-ib-1/2} (xy)^{-\a/2} \exp\left( -\frac{|x-y|^2}{4t}\right) \dtt} \\
=& \pi^{-1/2} {2^{2bi+1}} \Gamma\eee{ib+1/2} \cl (xy)^{-\a/2}|x-y|^{-2ib-1}.
}

Notice that $A_{1,1}$ is one of the terms from \eqref{kernel_eq}. Next, by \eqref{int-conv},
\alx{
    |R_2| &\leq C \cg \, |x-y|^{-1} (xy)^{-\a/2} \int_{\frac{|x-y|^2}{xy}}^\8 t^{1/2} e^{-t/4}\dtt \\
    &\leq Cxy(x+y)^{-\a-3},\\
	 |R_3| &\leq C \cl x^{-\a} |x-y|^{-1} \int_0^{\frac{|x-y|^2}{xy}} t^{1/2}  \dtt  \\
	& \simeq C \cl x^{-\a-1}	.
}

Now let us turn to study $A_2$. Denote $c_\a = 4^{-(\a-1)/2} \Gamma((\a+1)/4)^{-1}$. Then, by using \eqref{Bessel_loc},
\spx{
A_2 = &c_\a \int_{xy}^\8 t^{-ib-(\a+1)/2} \exp\left( -\frac{x^2+y^2}{4t}\right)\dtt +R_4 = A_{2,1}+R_4\\
}
where, by \eqref{int-conv},
\alx{|R_4| &= \abs{ \int_{xy}^\8 t^{-ib-1} (xy)^{-(\a-1)/2} \exp\left( -\frac{x^2+y^2}{4t}\right) \left(I_{(\a-1)/2}\left(\frac{xy}{2t}\right) - \Gamma\eee{\frac{\a+1}{4}}^{-1} \left(\frac{xy}{4t}\right)^{(\a-1)/2} \right)\dtt}\\
&\leq C xy \int_{xy}^\8 t^{-(\a+3)/2}  \exp\left( -\frac{x^2+y^2}{4t}\right)  \dtt \\
&\leq C  xy (x^2+y^2)^{-(\a+3)/2} \simeq C xy (x+y)^{-\a-3}.
}
Moreover,
\spx{
A_{2,1} +R_5 &= c_\a \int_0^\8 t^{-ib-(\a+1)/2} \exp\eee{-\frac{x^2+y^2}{4t}} \dtt \\
&=c_\a 4^{ib+(\a+1)/2} \Gamma\eee{ib+(\a+1)/2} \eee{x^2+y^2}^{-ib-(\a+1)/2} ,
}
where
\alx{
\abs{R_5} &= \abs{c_\a \int_0^{xy}t^{-ib-(\a+1)/2} \exp\left( -\frac{x^2+y^2}{4t}\right)\dtt}\\
&\leq C (x^2+y^2)^{-(\a+1)/2} \int_{\frac{x^2+y^2}{xy}}^\8 t^{(\a+1)/2} e^{-t/4} \dtt\\
&\leq C xy(x+y)^{-\a-3}.
}
}

\pr{[Proof of Theorem \ref{lower1} for $\a<0$]

Let {$|b|>1$} and $\e\in(0,10^{-1})$ (to be fixed later on).  Denote $I=[1, 1+\e]$ and $S=[1+3\e,2]$.
Put $f_\e(x) = \e^{-1} \chi_I(x) x^{-\a}$, so that $\norm{f_\e}_{L^1(X)}=1$. If $x\in S$, by Lemma \ref{kernel} and the triangle inequality,
\sp{ \label{neg_triangle}
 \abs{B^{ib} f_\e(x)} & \leq \abs{ c_2(b) } x^{-\a/2} |x-1|^{-1}  \\
 & + \abs{c_2(b)}\abs{\int_I \left( (xy)^{-\a/2}|x-y|^{-2ib-1} - x^{-\a/2}|x-1|^{-2ib-1} \right) f_\e(y) \, d\nu(y)} \\
 & + \abs{ c_1(b)}\abs{\int_I \left( x^2 + y^2 \right)^{-ib-(\a+1)/2} f_\e(y) \, d\nu(y)} \\
 & + \abs{c_3(b) }\abs{\int_I R_b(x,y) f_\e(y) \, d\nu(y)}\\
 & = \abs{ c_2(b) } x^{-\a/2} |x-1|^{-1} + \Lambda_1 + \Lambda_2 + \Lambda_3.
}
Observe that for $x\in S$ and $y\in I$ we have $|x-y|\simeq|x-1|$ and $x\simeq y\simeq 1$. By using the Mean Value Theorem for the function $y \mapsto y^{-\a/2} |x-y|^{-2ib-1}$,
\sp{ \label{neg_est}
\Lambda_1 & \leq C \abs{c_2(b)} \e^{-1} \int_{1}^{1+\e} |b||y-1| |x-1|^{-2} \, dy \leq C \e \abs{bc_2(b)} |x-1|^{-2},\\
\Lambda_2 & \leq C \abs{c_1(b)} \e^{-1} \int_{1}^{1+\e} (x^2+y^2)^{-(1+\a)/2} \, dy  \simeq C \abs{c_1(b)}, \\
\Lambda_3 & \leq C \abs{c_3(b)} \e^{-1} \int_{1}^{1+\e} xy(x+y)^{-\a-3} \, dy \simeq C \abs{c_3(b)}.
}
Fix $|b|\geq 1$ and $\la$ such that $\la > \max(\Lambda_2,\Lambda_3, |b c_2(b)| )$. Recall that $x^{-\a/2} \geq 1$ for $x\in S$, so that for $\e$ small enough
\sp{\label{la2}
\nu\set{x\in S: \abs{c_2(b)}x^{-\a/2}|x-1|^{-1} > 4\la} & \geq \nu\set{x\in S: \abs{c_2(b)}|x-1|^{-1} > 4\la} \\
& = \int_{1+3\e}^{1+\abs{c_2(b)}/(4\la)} x^{\a} \, dx\\
& \geq C  {\abs{c_2(b)}/(4\la)}}
and
\sp{\label{la1}
\nu \set{ x\in S: \abs{\Lambda_1} > \la }  & \leq \nu \set{ x\in S: C \e  \abs{bc_2(b)} |x-1|^{-2} > \la }\\
& \leq \int_{1+3\e}^{1 + C(\e \la^{-1}  \abs{bc_2(b)})^{1/2}} x^\a \, dx \\
& \leq C \left( \e \la^{-1} \abs{bc_2(b)}\right)^{1/2} \\
& \leq C \e^{1/2}.
}

Hence, using \eqref{neg_triangle}--\eqref{la1} and  \eqref{gamma_constants} we get
\spx{\label{neg_end}
\norm{B^{ib}f_\e}_{L^{1,\8}(X)} \geq & \la \nu \set{ x\in S: \abs{B^{ib}f_\e(x)} > \la  }  \geq \la \nu\set{x\in S: \abs{c_2(b)}x^{-\a/2}|x-1|^{-1} > 4\la }\\
&-  \la \nu \set{x\in S: \abs{\Lambda_1} > \la }
 -  \underbrace{\la \nu \set{x\in S: \abs{\Lambda_2} > \la } }_{=0}-   \underbrace{\la \nu \set{x\in S: \abs{\Lambda_3} > \la }}_{=0} \\
& {\geq C |c_2(b)| -
C\la \e^{1/2} \geq C |c_2(b)|} \simeq |b|^{1/2} = |b|^{d/2}.
}
}
}

{{Turning} to the case $\a>0$ we could also use Lemma \ref{kernel}. In this case, {the summand with $c_1(b)$ would play the first role. An alternative proof that we shall present here uses integral  representation of the modified Bessel function. The same will be used in the proof of } Theorem \ref{lower2}. It is known that for} $\a>0$
\eq{ \label{bint}
	I_{(\a-1)/2}(z) = \eee{\Gamma(\a/2)^{-1}\sqrt{\pi}}^{-1} \left( \frac{z}{2}\right)^{(\a-1)/2} \int_{-1}^1 e^{-zs}(1-s^2)^{\a/2-1} ds, \quad z>0,
}}
see \cite[Ch. 6]{Watson}. Therefore, for $\a>0$, using \eqref{Kb}, \eqref{Kt}, and \eqref{bint} we obtain
\sp{\label{kernel2}
 K_b(&x,y) = { (2\Gamma(-ib))^{-1}\int_0^\8 t^{-ib-1} (xy)^{-(\a-1)/2} I_{(\a-1)/2}\eee{\frac{xy}{2t}} \exp\left(-\frac{x^2+y^2}{4t} \right) \dtt} \\
&  {= \eee{2^\a \Gamma(-ib) \Gamma(\a/2)\sqrt{\pi}}^{-1}  \int_{-1}^1 \int_0^\8 t^{-ib-(\a+1)/2} \exp\left(-\frac{x^2+y^2+2xys}{4t} \right) \dtt \, (1-s^2)^{\frac{\a}{2}-1}  ds }\\
& {= \frac{2^{2ib+1} \Gamma\left( ib+ (\a+1)/2\right)}{\Gamma(-ib) \Gamma(\a/2)\sqrt{\pi}}  \int_{-1}^1 \left( x^2+y^2+2xys \right)^{-ib-(\a+1)/2}  \eee{1-s^2}^{\a/2-1} \, ds} \\
& = {C_\a} c_1(b) \int_{-1}^{1} \left( x^2+y^2+2xys \right)^{-ib-(\a+1)/2} \eee{1-s^2}^{\a/2-1} ds,
}
{where $C_\a = \pi^{-1/2}  \Gamma(\a/2)^{-1} \Gamma((\a+1)/4)$. }

\pr{[Proof of Theorem \ref{lower1} for $\a> 0$]
Let {$|b|>1$}, $\e\in(0,10^{-1})$, and $f_\e(x) = x^{-\a} \e^{-1} \chi_{[\e,2\e]}(x)$. Similarly as in \eqref{neg_triangle}, using \eqref{kernel2} and the Mean Value Theorem, for $x>3\e$ we have
 \sp{ \label{pos_triangle}
  \abs{B^{-ib}f_\e(x)}  \leq &  \abs{ \int_\e^{2\e} K_b(x,0)  f_\e(y) \, d\nu(y) } + \abs{ \int_\e^{2\e} \left( K_b(x,0) - K_b(x,y)  \right) f_\e(y) \, d\nu(y) }\\
  \leq & C|c_1(b)| \e^{-1} \left\{ \abs{\int_\e^{2\e}\int_{-1}^{1} \eee{1-s^2}^{\a/2-1} x^{-2ib-(\a+1)}\, ds \, dy
  }\right. \\
  & + \left.\int_\e^{2\e} \int_{-1}^{1} \eee{1-s^2}^{\a/2-1} \abs{ \left( x^2+y^2+2sxy \right)^{-ib-\frac{\a+1}{2}} - x^{-2ib-(\a+1)} } \, ds \, dy \right\}\\
    \leq & C|c_1(b)| \eee{ x^{-\a-1} +   \e^{-1} \int_{-1}^{1} \eee{1-s^2}^{\a/2-1} \int_\e^{2\e} |b|\abs{y^2+2sxy} x^{-\a-3} \, dy \, ds }  \\
 \leq  & C|c_1(b)| \eee{ x^{-\a-1} + \e \abs{b} x^{-\a-2}}.
  }

 Let us fix $|b|>1$ and $\la>|b c_1(b)|$. For all $\e$ small enough, we get
 \spx{
 \nu\set{x>3\e: C \abs{c_1(b)} x^{-\a-1} > 2\la } & = \int_{3\e}^{ C (\abs{c_1(b)}/\la)^{1/(1+\a)}} x^{\a} \, dx \geq  C \abs{c_1(b)}/\la,
 }
and
 \spx{
 \nu\set{x>3\e:  C \abs{bc_1(b)} \e x^{-\a-2} > \la}  \leq & \int_0^{ C \left(\abs{bc_1(b)} \e / \la\right)^{1/(\a+2)}} x^{\a} \, dx \\
  =  &C \left(\abs{bc_1(b)} \e/ \la\right)^{(\a+1)/(\a+2)} \\
 \leq  & C \e^{(1+\a)/(\a+2)}.
 }
 Therefore, by choosing a proper $\e$, we obtain

 \spx{ \label{pos_end}
 \norm{B^{ib}f_\e}_{L^{1,\8}(X)} \geq & \la \nu\set{ x\in X: \abs{B^{ib}f_\e(x)} > \la } \\
 \geq & \la \nu \set{ x > 3\e : C|c_1(b)| x^{-\a-1} > 2\la } - \la \nu \set{ x > 3\e :  C |c_1(b)| \e x^{-\a-2} > \la } \\
  \geq & C |c_1(b)| -  C \la \e^{(\a+1)/(\a+2)}  \geq  C |c_1(b)| \simeq |b|^{(\a+1)/2} = |b|^{d/2} .
 }
}

\pr{[Proof of Theorem \ref{lower2}]
Set $\a>0$, $p\in(1,2)$, $|b|>1$, $\e\in(0,10^{-1})$, $\de >1$, and a function ${f}\in{L^p(X)}$ such that $\supp f \subseteq (0,\e)$, and $f\geq 0$. Similarly as in \eqref{pos_triangle}, using \eqref{kernel2} and Corollary \ref{gamma2},
\spx{
\norm{\Bb^{ib}f}&_{L^p(X)}^p
\geq \int_\de ^\8 \abs{ \int_X \eee{K_b(x,0)- (K_b(x,0) - K_b(x,y)} f(y)\, d\nu(y)}^p\, d\nu(x) \\
\geq & C \norm{f}_{L^1(X)}^p \int_\de^\8 \abs{K_b(x,0)}^p  \, d\nu(x) - C \int_\de^\8 \abs{\int_X (K_b(x,y)-K_b(x,0)) f(y)\, d\nu(y)}^p \, d\nu(x) \\
 \geq & C  \norm{f}_{L^1(X)}^p|b|^{p(\a+1)/2} \eee{ \int_\delta^{\8} x^{-p(\a+1)+\a} \, dx
 - \int_\delta^{\8} \e^p |b|^p x^{-p(\a+2)+\a} \, dx} \\
  \geq & C \norm{f}_{L^1(X)}^p |b|^{p(\a+1)/2} \delta^{(\a+1)(1-p)} \eee{ 1 - \e^p |b|^{p} \delta^{-p}}.
}
Now we take $\delta=|b|$ and fix $\e$ small enough, independent of $b$, getting
\spx{\label{lpnorm}
	\norm{\Bb^{ib}f}_{L^p(X)} \geq C_{p} {|b|^{\frac{(\a+1)(2-p)}{2p}}} \norm{f}_{L^1(X)} \geq C_{p,\e} {|b|^{{\frac{d}{2}}\frac{2-p}{p}}} \norm{f}_{L^p(X)}.
}
}

\section{Appendix - Gamma function estimate}\label{app}

\lem{gamma}{
Let $a+bi\in \CC$. {For $a\geq 0$ fixed and all $|b|\geq 1$ we have}
	\eqx{
    |\Gamma(a+bi)| \simeq |b|^{a-1/2} \exp\eee{-\frac{\pi b}{2}}.
    }
}
The result above is known. It is a consequence of the Stirling's Formula, see \cite[Ch. 6]{Abramowitz_Stegun}. For the convenience of the reader we present a short proof.

\pr{
Using the reflection formula
\eq{ \label{reflection}
\Gamma(1-z)\Gamma(z) = \pi/\sin(\pi z),}
and the recursion identity
\eq{ \label{recursion}
z\Gamma(z)= \Gamma(z+1),}
we have that $|1-ib| \abs{\Gamma(ib)}^2 = \abs{\Gamma(ib)\Gamma(1-ib)} =  \abs{{\pi}/{\sin(\pi ib)}} \simeq \exp\eee{-\pi |b|}$ for $|b|\geq 1$. Thus,
\eq{\label{gamma0}
	\abs{\Gamma(ib)} \simeq |b|^{-1/2}\exp\eee{-\frac{\pi |b|}{2}}, \qquad |b|\geq 1
.}

Denote {$S = \set{z\in \CC \ : \ 1\leq \rm{Re}(z) \leq 2,  {\abs{\rm{Im}(z)}} \geq 1}$} and define a holomorphic function
\eqx{
F(z) = \Gamma(z)z^{-z+1/2}, \quad z\in S.}
Now, we claim that $|F(z)|\leq C$ if $z\in \partial S$. This is clear for $z=a\pm i$, $a\in[1,2]$. For $z=1+ib$, $|b|\geq 1$, we use \eqref{recursion} and \eqref{gamma0} getting
\spx{
	\abs{F(1+ib)} & = \abs{\Gamma\left(1+ib\right)} \abs{(1+ib)^{-1/2-ib}} = \abs{b}\abs{\Gamma\left(ib\right)} (1+b^2)^{-1/4}  e^{b\, \mathrm{arctg} b} \\
    & \leq C |b|^{1/2} e^{-\pi|b|/2} |b|^{-1/2} e^{b \, \rm{arctg}(b)}  \leq C.
}
Similarly we show boundedness of $F$ for $z=2+bi$, $|b|\geq 1$.

Observe that $\abs{F(z)} \leq |\Gamma(z)| |z|^{|-z+1/2|} \leq C e^{c|z|^2}
$ for $z\in S$. Hence, applying the Phragm\'{e}n-Lindel\"{o}f principle, we obtain that $|F(z)| \leq C$ for $z \in S$. Therefore, for a fixed $a\in[1,2]$ and $|b|\geq 1$ we have
\sp{\label{up_up}
\abs{\Gamma(a+bi)} \leq C \abs{(a+
bi)^{a-1/2+bi}} = \eee{a^2+b^2}^{(2a-1)/4} \cdot e^{-b \, \mathrm{arctg}(b/a)}
\simeq C|b|^{a-1/2} e^{-\pi |b|/2}.
}
This is the desired estimate from above for $a\in[1,2]$. We extend this for all $a\in[0,\8)$ by using \eqref{recursion}. Then, by \eqref{reflection}, we get estimate from below for $a\in[0,1]$, and extend this for $a\in [0,\8)$ using \eqref{recursion} once more.
}

\cor{gamma2}{
For fixed $a_1,a_2 \geq 0$ and $|b|\geq 1$ we have
	\eqx{
    \abs{\frac{\Gamma(a_1+bi)}{\Gamma(a_2+bi)}}\simeq |b|^{a_1-a_2}.
    }
}

{\bf Acknowledgments:} The authors would like to thank Jacek Dziuba\'nski, Alessio Martini, Adam Nowak, B\l a\.zej Wr\'obel, and the referees for their helpful comments and suggestions.

\bibliographystyle{amsplain}        

\def\cprime{$'$}
\providecommand{\bysame}{\leavevmode\hbox to3em{\hrulefill}\thinspace}
\providecommand{\MR}{\relax\ifhmode\unskip\space\fi MR }
\providecommand{\MRhref}[2]{%
  \href{http://www.ams.org/mathscinet-getitem?mr=#1}{#2}
}
\providecommand{\href}[2]{#2}

\end{document}